\documentclass{article}
\def\BibTeX{{\rm B\kern-.05em{\sc i\kern-.025em b}\kern-.08em
    T\kern-.1667em\lower.7ex\hbox{E}\kern-.125emX}}

\usepackage{amssymb,amsfonts,amsthm,amsmath}
\usepackage{fontawesome5}
\usepackage[colorlinks,allcolors=BrickRed]{hyperref}
\usepackage{textcomp}
\usepackage[dvipsnames]{xcolor}
\usepackage{algorithmic}
\usepackage{algorithm}
\usepackage{fancyvrb}
\usepackage{dsfont}

\usepackage{enumitem}
\usepackage{nicefrac}
\usepackage{schemabloc}
\usepackage{circuitikz}
\usepackage{comment}

\usetikzlibrary{circuits}

\newcommand{\A}{\mathbb{A}}

\renewcommand{\L}{\mathbb{L}}
\newcommand{\M}{\mathbb{M}}
\newcommand{\N}{\mathbb{N}}

\newcommand{\R}{\mathbb{R}}

\newcommand{\X}{\mathbb{X}}

\newcommand{\bbGamma}{{\mathpalette\makebbGamma\relax}}
\newcommand{\makebbGamma}[2]{%
  \raisebox{\depth}{\scalebox{1}[-1]{$\mathsurround=0pt#1\mathbb{L}$}}%
}

\newcommand{\bfA}{\mathbf{A}}

\newcommand{\bfx}{\mathbf{x}}

\newcommand{\rmd}{\mathrm{d}}

\newcommand{\rmt}{\mathrm{t}}

\newcommand{\rmx}{\mathrm{x}}

\newcommand{\cL}{\mathcal{L}}

\newcommand{\cQ}{\mathcal{Q}}
\newcommand{\cR}{\mathcal{R}}

\newcommand{\vx}{\boldsymbol{x}}
\newcommand{\vy}{\boldsymbol{y}}

\newcommand{\vf}{\boldsymbol{f}}
\newcommand{\vg}{\boldsymbol{g}}
\newcommand{\vh}{\boldsymbol{h}}

\newcommand{\vO}{\boldsymbol{0}}

\newcommand{\vi}{\boldsymbol{i}}
\newcommand{\vj}{\boldsymbol{j}}

\newcommand{\vu}{\boldsymbol{u}}
\newcommand{\vv}{\boldsymbol{v}}

\newcommand{\vphi}{\boldsymbol\phi}

\newcommand{\bTheta}{\boldsymbol\Theta}

\newcommand{\vpartial}{\boldsymbol\partial}

\newcommand{\vol}{\mathrm{vol}}
\newcommand{\ind}{\mathds{1}}

\newcommand{\x}{\times}

\newcommand{\dist}{\mathrm{dist}}

\newtheorem{prop}{Proposition}
\theoremstyle{definition}
\newtheorem{dfn}{Definition}
\newtheorem{asm}{Assumption}
\theoremstyle{remark}
\newtheorem{rem}{Remark}

\usepackage{geometry}
\usepackage[nolist]{acronym}

\begin{document}

\title{SOStab: a Matlab Toolbox for Transient Stability Analysis}

\author{St\'ephane Drobot\textsuperscript{1}, Matteo Tacchi\textsuperscript{2}\thanks{Corresponding author. {\scriptsize \faEnvelope } \texttt{matteo.tacchi@gipsa-lab.fr.} \faPhoneSquare* +33 4768 26235. \faMapMarker* 11 rue des Mathématiques, Grenoble Campus BP46, 38402 Saint-Martin-d’Hères Cedex, France.} {}, Carmen Cardozo\textsuperscript{1} and Colin N. Jones\textsuperscript{3}
}

\footnotetext[1]{R\&{}D division, Réseau de Transport d'Électricité, La Défense, France}

\footnotetext[2]{Univ. Grenoble Alpes, CNRS, Grenoble INP (Institute of Engineering Univ. Grenoble Alpes), GIPSA-lab, Grenoble, France.}

\footnotetext[3]{Autoamtic Control Laboratory, EPFL, Lausanne, Switzerland. 
}

\begin{acronym}
\acro{roa}[RoA]{Region of Attraction}
\acro{sos}[SoS]{Sums-of-Squares}
\acro{sosp}[SoSP]{SoS Programming}
\acro{eeac}[EEAC]{Extended Equal Area Criterion}
\acro{lmi}[LMI]{Linear Matrix Inequality}
\acro{sdp}[SDP]{Semidefinite Programming}
\end{acronym}

\maketitle

\begin{abstract}
This paper presents a new Matlab toolbox, aimed at facilitating the use of polynomial optimization for stability analysis of nonlinear systems. In the past decade several decisive contributions made it possible to recast this type of problems as convex optimization ones that are tractable in modest dimensions. However, available software requires their user to be fluent in Sum-of-Squares programming, preventing them from being more widely explored by practitioners. To address this issue, SOStab entirely automates the writing and solving of optimization problems, and directly outputs relevant data for the user, while requiring minimal input. In particular, no specific knowledge of optimization is needed for implementation. The toolbox allows a user to obtain outer and inner approximates of the \ac{roa} of the operating point of different grid connected devices such as synchronous machines and power converters.
\end{abstract}

\begin{center}
\textbf{Keywords}
\end{center} \small
AC power systems; Software; Sum-of-squares programming; Transient stability analysis. \normalsize

\begin{center}
\textbf{Acknowledgement}
\end{center} \small
This work was supported by the Swiss National Science Foundation under the NCCR Automation project, grant agreement 51NF40\underline{ }180545, and the French company RTE, under the RTE-EPFL partnership n$^\circ$2022-0225. \normalsize

\section{Introduction}

{A \textit{\ac{roa}} characterizes a set of initial conditions within the state space of a dynamical system for which trajectories demonstrate stability properties, such as hitting a target (finite time horizon \ac{roa}~\cite{henrionOuterROA,kordaInnerROA}) or converging towards a secure zone of the state space (infinite time horizon \ac{roa}~\cite{Wloszek,kordaMPI,oustryMPI}). Practical computation of approximations for both types of \ac{roa} is achievable by solving a \ac{sosp} problem within the framework known as Lasserre's Moment-\ac{sos} hierarchy, which involves resolving \ac{sdp} problems. While accommodating nonlinear behaviors, this approach comes with convergence guarantees~\cite{LasserreBook,MomentSOS,TacchiCV}, providing an infinite set of \textit{certified} stable operating conditions through a singular computation.} 

{Consequently, \ac{roa} has captured the interest of the power systems and control community. Notably, the potential application of \ac{sosp} to the transient stability problem have been explored over the past decade~\cite{anghel,Marinescu,Izumi,joszTSA,oustry}. Typically performed through contingency analyses using time-domain simulations, transient stability assessment is becoming increasingly challenging as the system operates near its limits in an environment of growing variability and uncertainty.} 

{In this scenario, having a means to rapidly eliminate situations with guaranteed stability proves invaluable. To serve this purpose, direct methods relying on simplified models, including those based on energy functions or on the \ac{eeac}~\cite{eeac}, are currently being revisited. In this context, \ac{roa}-based approaches emerge as a natural alternative, with scalability issues being diligently addressed by the applied mathematics community (see~\cite{tacchiThese} for details and~\cite{tacchi,subotic,wang,josz,franc} for existing solutions).}

{In this work, we address a set of more practical barriers that currently limit further development in this direction, such as:}
\begin{itemize}
    \item {The lack of a user-friendly interface for a practical implementation and application of the Moment-\ac{sos} hierarchy.}
    \item {The conditioning of the underlying \ac{lmi}, which depends on the specific formulation of the problem.}
\end{itemize}

More precisely, the existing frameworks~\cite{yalmip,GloptiPoly,SOSTOOLS} require users to manually write and solve \ac{sosp} problems to obtain the \ac{roa} approximation. In contrast, the proposed SOStab Matlab toolbox fully automates the \ac{sosp} aspect, eliminating the need for users to possess knowledge of the Moment-\ac{sos} hierarchy. The toolbox operates with minimal input requirements, namely dynamics, state constraints, equilibrium point, time horizon, target set, and a complexity parameter $d$. It outputs the stability certificate describing the \ac{roa} approximation and provides graphical representation in selected state coordinates. SOStab has been developed and made \href{https://github.com/droste89/SOStab}{publicly available} which allows users to compute \ac{roa} of non-linear dynamical systems.

{To the best of the authors' knowledge, the only existing toolbox for \ac{roa} approximation in a finite time horizon setting is \href{https://github.com/wangjie212/SparseDynamicSystem}{SparseDynamicSystem}~\cite{wang}, coded in the Julia language. A notable distinction lies in the fact that the Julia toolbox exclusively supports polynomial dynamics, whereas SOStab is built on the Matlab codes supporting~\cite{joszTSA}, specifically tailored for AC power systems that include phase variables and trigonometric dynamics.}

The article is organized as follows: first Section~\ref{sec:roa}
provides theoretical background on the definition of \ac{roa} and the calculation of inner and outer approximations, introducing the methods behind the tool and their relevance to the power system transient stability problem. 
Basics about the \ac{sosp} framework are included in Appendix~\ref{sec:hierarchy}. Then, Section~\ref{sec:modeling} presents the models of the two case studies considered in this work for illustrative purposes. Section~\ref{sec:tuto} includes an overview of the SOStab toolbox and offers guidance on its usage, while simulation results are discussed in Section~\ref{sec:results}. Finally, Section~\ref{sec:conclusion} summarizes the main contributions of this work and elucidates on potential improvements.  

\section{Methodological Framework}
\label{sec:roa} 

Consider a generic differential system
\begin{subequations} \label{eq:syst}
\begin{equation} \label{eq:ODE}
    \dot{\vx}(t) = \vf(\vx(t))
\end{equation}
with polynomial dynamics $\vf \in \R[\vx]$, equilibrium $\vx^\star \in \R^n$. Define a threshold $\Delta\vx \in (0,\infty)^n$ and state constraint
\begin{equation} \label{eq:state_con}
    \vx(t) \in \X := [\vx^\star \pm \Delta\vx],
\end{equation}
\end{subequations}
\small $[\vx^\star \pm \Delta\vx] := [x^\star_1 - \Delta x_1, x^\star_1 + \Delta x_1]\x\ldots\x[x^\star_n -\Delta x_n, x^\star_n + \Delta x_n]$.\normalsize

\begin{dfn}[{\ac{roa}~\cite{henrionOuterROA,kordaInnerROA}}] \label{dfn:RoA}
Given a time horizon $T > 0$ and a closed target set $\M \subset \X$, the Region of Attraction $\cR_T^\M$ of $\M$ in time $T$ is defined as {the set of all initial conditions of trajectories that hit the target $\M$ at time $T$, i.e.:}
\begin{equation} \label{eq:RoA}
    \cR_T^\M := \left\{\vx(0) \in \R^n \; : \begin{array}{c}
        \forall t \in [0,T],~\eqref{eq:syst} \text{ holds} \\
        \vx(T) \in \M
    \end{array} \right\}.
\end{equation}
\end{dfn}

{The transient stability problem can be reformulated as the search for a \ac{roa} in the state space of the post-fault system. Determining the duration a fault can last before a synchronous generator loses stability is equivalent to identifying how \textit{far} the trajectories can diverge before they can no longer return to an acceptable operating point. Here, the vector field $\vf$ models the nonlinear system dynamics, and the target $\M$ represents a positively invariant vicinity of a stable equilibrium point $\vx^\star$, both in the post-disturbance state. Consequently, the term \textit{initial condition} refers to the system state at the moment of fault clearance, and $T$ characterizes, in some manner, the recovery time. Finally, state constraints in Eq.~\eqref{eq:state_con} set bounds on state variables accounting for system limits.

Notice that here the dynamcis $\vf$ are assumed to be polynomial: for non-polynomial dynamics (e.g. phase-related sines and cosines), a preprocessing (Taylor expansion or variables change) is required to enforce polynomial structure.}

SOStab {builds on \cite{henrionOuterROA,kordaInnerROA} and} is specifically designed to {approximate} such an $\cR_T^\M$ for ellipsoids $\M$ described by
\begin{equation} \label{eq:target}
    \M := \{\vx \in \R^n \; : \; \|\bfA(\vx-\vx^\star)\| \leq \varepsilon \}
\end{equation}
where $\bfA \in \R^{n\x n}$ is a reshaping matrix such that $\det(\bfA) = 1$, and $\varepsilon > 0$ is an error tolerance.

From a computational viewpoint, considering finite time horizon \ac{roa} improves the properties of the resulting \ac{sos} programs: they are convex, while their infinite time horizon counterparts are bilinear (and thus nonconvex), which has important consequences on the convergence properties of the corresponding algorithms (see e.g.~\cite{Izumi}).
In practice, this is a reasonable assumption as operational points of physical system, such as power systems, are constantly moving, making the results of the analysis relevant for a limited time window. 
In exchange for the finite time horizon, it is necessary to consider target sets that are not reduced to a point, defined by parameters $\bfA$ and $\varepsilon$; in SOStab, the default value for $\bfA$ is the identity matrix, but sometimes (e.g. when variables evolve on different time scales), it is useful to make a different choice. { For instance, in a two-dimensional singular perturbation framework $|\dot{x}_1| \sim a^2|\dot{x}_2|$, one would set 
$$\bfA = \begin{pmatrix}
    \nicefrac{1}{a} & 0 \\ 0 & a
\end{pmatrix}$$
i.e. take an ellipse of equation $\nicefrac{(x_1-x_1^\star)^2}{a^2} + a^2 (x_2-x_2^\star)^2 \leq \varepsilon^2$ for $\M$, which is justified by the physical reasoning that, for a trajectory starting at $(x_1,x_2)$ and converging to $(x_1^\star,x_2^\star)$ in infinite time, it holds $$x_1^\star - x_1 = \int_0^\infty \dot{x}_1(t) \; dt \sim a^2 \int_0^\infty \dot{x}_2(t) \; dt = a^2(x_2^\star - x_2),$$
and those two terms should have equal contribution in the description of the target $\M$, otherwise the resulting \ac{roa} will be skewed. Also, from Eq.~\eqref{eq:target} it holds that $\vol(\M) \propto \varepsilon^n$, where $\mathrm{vol}$ denotes the $n$-dimensional volume of a set, so that the choice of the tolerance parameter $\varepsilon$ is directly related to the desired size of the target $\M$.}

SOStab takes as input the problem data $(\vf,\Delta\vx,T,\bfA,\varepsilon)$, as well as a parameter $d \in 2\N$ {which sets the accuracy and complexity of the approximation}, and outputs the following objects for outer and inner \ac{roa} estimation:
\begin{subequations} \label{eq:output}
    \begin{align}
        \vy_d^{out} & = \begin{pmatrix} \lambda_d^{out} & v_d^{out} & w_d^{out}  \end{pmatrix} \label{eq:outputout} \\
        \vy_d^{in} & = \begin{pmatrix} \lambda_d^{in} & v_d^{in} & w_d^{in} \end{pmatrix}, \label{eq:outputin}
    \end{align}
\end{subequations}
with $\lambda_d^{in / out} \geq 0$, $v_d^{in / out} \in \R_d[t,\vx]$ degree $d$ polynomials in $(t,\vx)$ and $w_d^{in / out} \in \R_d[\vx]$ degree $d$ polynomials in $\vx$, such that the following inclusion guarantees hold:
\begin{subequations} \label{eq:guarantee}
    \begin{align}
        \cR_T^\M & \subset \{\vx \in \R^n \; : \; v_d^{out}(0,\vx) \geq 0\} =: \cR_d^{out} \label{eq:outer} \\
        \cR_T^\M & \supset \{\vx \in \R^n \; : \; v_d^{in}(0,\vx) < 0 \} =: \cR_d^{in}.  \label{eq:inner}
    \end{align}
\end{subequations}
In words, {$v_d^{out}(0,\vx) \geq 0$ (resp. $v_d^{in}(0,\vx)<0$) means that $\vx$ belongs to} an outer (resp. inner) approximation $\cR_d^{out}$ (resp. $\cR_d^{in}$) of $\cR_T^\M$. 
Moreover, 
{
$\lambda_d$ is an upper bound of the approximation error, and $w_d$ is often used to check numerical results: it should separate $\X$ into a zone where its values are close to $0$ and one where they are above $1$; one of those two zones missing means that the solver failed to detect the \ac{roa}}. Eventually, the framework comes with precision guarantees:
\begin{equation}
    \vol\left(\cR_d^{in / out}\right) \underset{d\to \infty}{\longrightarrow} \vol\left(\cR_T^\M\right). \label{eq:precio}
\end{equation}
 In words, when the modelling parameter $d$ goes to infinity, the volume of the error made by the approximation schemes vanishes (see Appendix~\ref{sec:hierarchy} or~\cite{henrionOuterROA,kordaInnerROA} for further details { and Appendix~\ref{sec:pol} for an elementary example}).

{
In the transient stability problem, inner approximations are more relevant as they provide stability certificates for all states within, even if it means excluding some post-fault conditions leading to stable trajectories. Outer approximations ensure the inclusion of all stable states but lack conservatism guarantees. The gap between the two may inform about the approximation accuracy with respect to the exact \ac{roa}.}

\section{Models for case study}
\label{sec:modeling}

In this section we present the model of two different test cases, proposed in the literature, to illustrate the tool capabilities, flexibility and performance. The first one is the well known synchronisation loop of classic grid-following converters: the phase locked loop (PLL~\cite{polypll,sync,frontiers}). The second {consists of a synchronous machine described by a third order model, to which a governor and automatic voltage regulator (AVR) representations are added. The grid side is represented by an infinite bus~\cite{Marinescu,oustry}.}

\subsection{Phase Locked Loop 2$^{\text{nd}}$ order model}

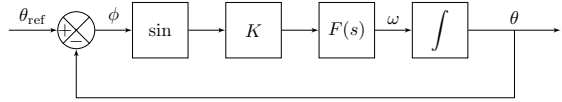
\begin{figure}[htbp]
    \centering
    \scalebox{0.7}{
    \begin{tikzpicture}
        \sbEntree{E}
        \sbComp{a}{E}
        \sbBloc{b}{$\sin$}{a}
                \sbRelier[$\theta_{\rm ref}$]{E}{a}
        \sbBlocL{c}{$K$}{b}
                \sbRelier[$\phi$]{a}{b}
        \sbBlocL{d}{$F(s)$}{c}
        \sbBloc{e}{$\displaystyle\int$}{d}
                \sbRelier[$\omega$]{d}{e}
        \sbSortie[5]{S}{e}
                \sbRelier[$\theta$]{e}{S}
        \sbRenvoi{e-S}{a}{}
    \end{tikzpicture}
    }
    \caption{Block scheme of a PLL system.}
    \label{fig:pll}
\end{figure}

Fig.~\ref{fig:pll} shows a generic PLL block diagram~\cite{polypll}. Here an angular state variable $\theta$ is required to match a reference $\theta_{\rm ref}$; to that end, the system computes the sine of the phase difference $\phi$, multiplies it by some gain $K$ and takes it through a low-pass filter with transfer function $F(s) = \frac{1 + \tau_2 s}{\tau_1 s}$, resulting in the following differential system:
\begin{equation} \label{eq:pll}
    \begin{pmatrix} \dot\phi \\ \dot\omega \end{pmatrix} = \begin{pmatrix}
        \omega \\
        -K \frac{\tau_2}{\tau_1}\cos(\phi) \omega - \frac{K}{\tau_1} \sin(\phi) 
    \end{pmatrix}.
\end{equation}

Following~\cite{polypll}, in the PLL setting the time constants $\tau_1$ and $\tau_2$ are functions of the gain $K$, a natural frequency $\omega_n$ and a damping ratio $\zeta$:
\begin{equation}
    \tau_1 = \frac{K}{\omega_n^2} \qquad \tau_2 = \frac{2\zeta}{\omega_n}
\end{equation}

{
In our case study, we will use SOStab to compute the \ac{roa} of the PLL system described by Eq.~\eqref{eq:pll} with parameter values as described in Tab.~\ref{tab:pllparam}, and where sine and cosine have been replaced by their degree $10$ Taylor expansions (denoted \texttt{si} and \texttt{co} respectively), so that the dynamics are polynomial.

\begin{table}[htbp]
    \centering 
    \begin{tabular}{|c|c|c|c|}
        \hline
        Parameter & $K$ & $\omega_n$ & $\zeta$ \\
        Value & $1$ s$^{-1}$ & $10.813$ s$^{-1}$ & $1.3303$ \\
        \hline
    \end{tabular}
    \vspace{1em}
    \caption{Parameters for phase locked loop.}
    \label{tab:pllparam}
\end{table}
}

\subsection{Single machine - infinite bus with regulations}

Following~\cite{Marinescu}, Fig.~\ref{fig:nomeq} represents a synchronous machine connected to an infinite bus with voltage $\vv_s$ through a power transmission line with impedance $Z_t = R_t + j \, X_t$. The physics of the synchronous machine are modelled by the swing equation as well as the electromotive force dynamics:
\begin{subequations} \label{eq:SMIB}
    \begin{align}
        \dot{\theta} & = \omega - \omega_s \label{eq:phase} \\
        2H \dot{\omega} & = P_m - (v_d i_d + v_q i_q + ri_d^2 + ri_q^2) \label{eq:freq} \\
        T'_{d_0} \dot{e}'_q & = -e'_q - (x_d-x'_d)i_d + E_{fd}
    \end{align}
where $\theta$ is the machine angle in a rotating frame synchronized with the grid (hence the grid frequency $\omega_s$ in \eqref{eq:phase}) and the current $\vi$ and voltage $\vv$ are described in the rotating frame by the Kirchhoff laws:
\begin{align*}
    & i_q = \frac{(X+x'_d)\|\vv_s\|\sin\theta - (R+r)(\|\vv_s\|\cos\theta - e'_q)}{(R+r)^2+(X+x'_d)(X+x_q)} \\ 
    & i_d = \frac{X+x_q}{R+r}i_q - \frac{\|\vv_s\|\sin\theta}{R+r} \\ 
    & v_d = x_qi_q - ri_d \\ 
    & v_q = Ri_q + X i_d + \|\vv_s\|\cos\theta 
\end{align*}

\begin{figure}[htbp]
\centering 
\begin{circuitikz}[scale = 0.7, transform shape] \draw
   node[circ] (A) at (0,2) {\Large{$^A$}}
   node[ground] (Ta) at (0,0) {}
   node[oscillator] (G) at (-1,2){}
   (G) to[short] (A)
   (1,2) to[R, l=$Z_t$, i>=$\vi$] (7,2){}
   node[circ] (B) at (8,2) {\Large{$^B$}}
   (B) to[short] (9,2)
   node[ground] (Tb) at (8,0) {}
   ;
   \draw[-triangle 45] (7.8,0.3) -- (0.2,0.3) node[below, midway] {$\vu_t$};
   \draw[-triangle 45] (0,0) -- (0,1.9) node[left, midway] {$\vv$};
   \draw[-triangle 45] (8,0) -- (8,1.9) node[right, midway] {$\vv_s$};
   \draw (0.05,2) -- (1,2);
   \draw (7,2) -- (7.95,2);
   \draw (9.5,2) [thick]circle (0.5cm) node {\Large{$\infty$}};
; \end{circuitikz}
\caption{A synchronous machine connected to an infinite bus.}
\label{fig:nomeq}
\end{figure}
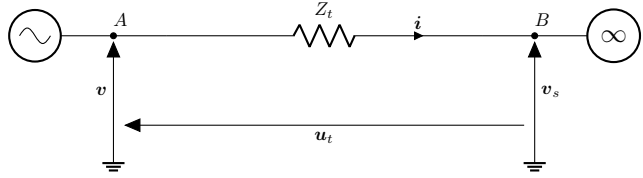

Moreover, the model includes quadratic AVR dynamics
\begin{equation} \label{eq:avr}
    T_a \dot{E}_{fd} = -E_{fd} + K_a(V^2-\|\vv\|^2)
\end{equation}
as well as turbine governor equation
\begin{equation} \label{eq:gov}
    T_g \dot{P}_m = -P_m + P + K_g(\omega_s-\omega)
\end{equation}
\end{subequations}
where $V$ and $P$ are given set points for the voltage magnitude $\|\vv\|$ and mechanical power $P_m$. In our case study, we will use SOStab to compute the \ac{roa} of the 5D SMIB model  given by Eq.~\eqref{eq:SMIB} with parameter values as described in Tab. \ref{tab:smibparam}, along with the variable change 
\begin{equation} \label{eq:variable}
    \vx = (\sin \theta, \cos \theta, \omega, e'_q, E_{fd}, P_m)
\end{equation}
so that the dynamics $\dot{\vx} = \vf(\vx)$ are polynomial.

\begin{table}[htbp]
    \centering 
    \begin{tabular}{|c|c|c|c|c|c|}
        \hline
        Parameter & $T'_d$ & $x_d$ & $x'_d$ & $x_q$ & $r$ \\
        Value & $9.67$ s & $2.38$ pu & $0.336$ pu & $1.21$ pu & $0.002$ pu \\
        \hline \hline
        Parameter & $H$ & $\omega_s$ & $R$ & $X$ & $V = \|\vv_s\|$ \\
        Value & $3$ s & $1$ pu & $0.01$ pu & $1.185$ pu & $1$ pu \\
        \hline \hline
        Parameter & $T_a$ & $K_a$ & $T_g$ & $K_g$ & $P$ \\
        Value & $1$ s & $70$ pu & $0.4$ pu & $0.5$ pu & $0.7$ pu \\
        \hline
    \end{tabular}
    \vspace{1em}
    \caption{Parameters for single machine - infinite bus.}
    \label{tab:smibparam}
\end{table}

\section{Computing RoA with SOStab}
\label{sec:tuto}

In this section we showcase a step by step procedure to assess the stability of our test cases with SOStab.

\color{black}

\subsection{Installation} \label{sec:instal}
SOStab is a freeware subject to the General Public Licence (GPL) policy available for Matlab. It can be downloaded at:
\begin{center}
\texttt{\url{https://github.com/droste89/SOStab}}
\end{center}

SOStab requires YALMIP~\cite{yalmip}, as well as a semidefinite solver. Mosek~\cite{mosek} is used by default, but it can be replaced by any other solver, provided they are installed and interfaced through YALMIP.

\subsection{Overview} \label{sec:overview}

SOStab is designed to compute a classifier, denoted as $v(0,\cdot)$, which can be evaluated at any state $\vx$ of the considered system. The sign of $v(0,\vx)$ indicates whether the state belongs or not to the selected \ac{roa} approximation, thereby characterizing the stability of the trajectory starting at $\vx$. \color{black}
The toolbox requires the following minimal input:
    \begin{align*}
        &\dot{\vx} = \vf(\vx), \quad \vf\in\R[\bfx]^n, \quad T>0 \\
        &\vx^\star \in \{\vx \in \R^n \; : \; \vf(\vx) = \vO\}\\
        & \Delta\vx \in (0,\infty)^n, \quad d \in 2\N, \quad \varepsilon > 0,
    \end{align*}
    {where $\vf \in \R[\bfx]$ means that the system dynamics should have a polynomial form (up to Taylor expansion or variable change).}
Then, it identifies the problem to be solved as: ``compute approximations $\cR_d$ of $\cR_T^\M$ with dynamics $\vf$, admissible set $[\vx^\star\pm\Delta\vx]$ and target set $\M = \{\vx \in \R^n \; : \; \|\vx - \vx^\star \| \leq \epsilon \}$ ".

SOStab also admits two optional input arguments:
\begin{itemize}
    \item a shape matrix $\bfA \in \R^{n\x n}$ such that $\det(\bfA)=1$ reshapes the target set as $\M = \{\vx \in \R^n \; : \; \|\bfA(\vx - \vx^\star) \| \leq \epsilon \}$; 
    \item in case the system at hand involves trigonometric functions of phase variables $\theta_1,\ldots,\theta_N$, it is also possible to specify a phase index matrix $\bTheta \in \R^{N\x 2}$ whose first (resp. second) column consists of the indices of the sines (resp. cosines) of the $\theta_i$ in the recasted variable $\vx$.
\end{itemize}

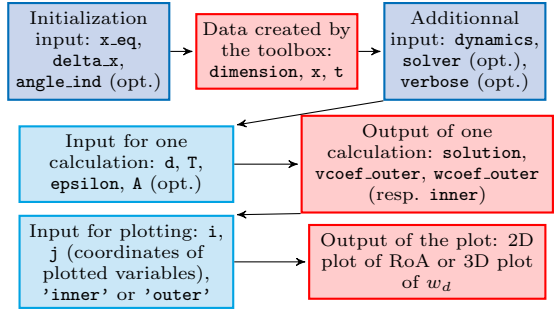
\begin{figure}[htbp]
    \centering
\begin{tikzpicture}
\tikzstyle{mybox} = [thick, rectangle, inner sep=2pt, inner ysep=3pt, font=\scriptsize]
\node [mybox, draw=NavyBlue, fill=NavyBlue!20](boxInput) at (-0.5,0){\parbox{2cm}{\centering Initialization input: \texttt{x\_eq}, \texttt{delta\_x}, \texttt{angle\_ind} (opt.)}};
\node [mybox, draw=red, fill=red!20](boxTool) at (2,0){\parbox{2cm}{\centering Data created by the toolbox: \texttt{dimension}, \texttt{x}, \texttt{t}}};
\node [mybox, draw=NavyBlue, fill=NavyBlue!20](boxAdd) at (4.5,0){\parbox{2cm}{\centering Additionnal input: \texttt{dynamics}, \texttt{solver} (opt.), \texttt{verbose} (opt.)}};
\node [mybox, draw=Cerulean, fill=Cerulean!20](boxInput2) at (0,-1.5){\parbox{2.7cm}{\centering Input for one calculation: \texttt{d}, \texttt{T}, \texttt{epsilon}, \texttt{A} (opt.)}};
\node [mybox, draw=red, fill=red!20](boxOutput1) at (4,-1.5){\parbox{3.2cm}{\centering Output of one calculation: \texttt{solution}, \texttt{vcoef\_outer}, \texttt{wcoef\_outer} (resp. \texttt{inner})}};
\node [mybox, draw=Cerulean, fill=Cerulean!20](boxInput3) at (0,-2.8){\parbox{2.7cm}{\centering Input for plotting: \texttt{i}, \texttt{j} (coordinates of plotted variables), \texttt{'inner'} or \texttt{'outer'}}};
\node [mybox, draw=red, fill=red!20](boxOutput2) at (4,-2.8){\parbox{3cm}{\centering Output of the plot: 2D plot of \ac{roa} or 3D plot of $w_d$}};
\draw[color=black,->,>=stealth',shorten >=1pt] (boxInput.east) -- (boxTool.west);
\draw[color=black,->,>=stealth',shorten >=1pt] (boxTool.east) -- (boxAdd.west);
\draw[color=black,->,>=stealth',shorten >=1pt,xshift=2pt,yshift=2pt](boxAdd.south west) -- (boxInput2.north east);
\draw[color=black,->,>=stealth',shorten >=1pt] (boxInput2.east) -- (boxOutput1.west);
\draw[color=black,->,>=stealth',shorten >=1pt] (boxOutput1.south west) -- (boxInput3.north east);
\draw[color=black,->,>=stealth',shorten >=1pt] (boxInput3.east) -- (boxOutput2.west);
\end{tikzpicture}
\caption{Flowchart of the toolbox workings}
\label{fig:flowchart}
\end{figure}

The use of the toolbox consists of four steps (see Fig.~\ref{fig:flowchart}):
\begin{enumerate}
    \item The initialization is performed from input $(\vx^\star,\Delta\vx)$ (and optional input $\bTheta$); it defines the admissible set $[\vx^\star\pm\Delta\vx]$ and identifies the dimension and variables of the system.
    \item The user then inputs the dynamics $\vf$ of the system and adjusts optional settings of the toolbox.
    \item The \ac{roa} approximation itself is performed from input $(d,T,\varepsilon)$ (with optional input $\bfA$); it outputs the \ac{sdp} solutions $\vy_d
    = (\lambda_d, v_d, w_d)$.
    \item Graphic representations of the solutions can eventually be plotted; the choice of the plots abscissa and ordinate and plotting options are up to the user.
\end{enumerate}

{We are now going to demonstrate the implementation of those four steps on the PLL and SMIB test cases.

\begin{rem} In the current version of the toolbox, when using an optional argument, one should also specify all optional arguments that appear before in the method call.
\end{rem}

\subsection{Initialization and dynamics specification} \label{sec:init}

The first initialization step consists in choosing the threshold $\Delta\vx$ defining the admissible set $\X$; For the PLL, we follow the reference~\cite{polypll} and ask $\Delta\phi =\pi$ rad and $\Delta\omega = 20\pi$ rad/s (the equilibrium is trivially $\vx^\star = \vO$). Similarly, we identify an equilibrium point $\vx^\star$ and set a threshold $\Delta \vx$ for the SMIB model (with variable $\vx$ described as in Eq.~\eqref{eq:variable}).

\begin{Verbatim}[frame=single, fontsize=\small]
  PLLeq = [0;0];   DeltaPLL = [pi; 20*pi];
  SMIBeq = [sin(1.539); cos(1.539); ...
              1; 1.070; 2.459; 0.7];
  DeltaSMIB = [1; 1; 1; 1; 20; 4];
\end{Verbatim}

Remember that the PLL dynamics are approximated by their Taylor expansion, while we study the exact SMIB model through a variable change. Hence, for the SMIB model, we need to specify which coordinates of $\vx$ are actually trigonometric functions of the original variable, by adding the input:

\begin{Verbatim}[frame=single, fontsize=\small]
  angle_ind = [1,2];
\end{Verbatim}

Here \texttt{angle\_ind} $=\bTheta$ is the phase index matrix with only one line as the SMIB model features only one phase angle $\theta$; the first (resp. second) column gives the position of $\sin\theta$ (resp. $\cos\theta$) in the recasted variable $$\vx = (\sin\theta, \cos\theta, \omega, e'_q, E_{fd}, P_m).$$

Next, for each of our two test cases we create an instance of SOStab with the following commands:

\begin{Verbatim}[frame=single, fontsize=\small]
  PLL = SOStab(PLLeq, DeltaPLL);
  SMIB = SOStab(SMIBeq, DeltaSMIB, Z);
\end{Verbatim}
}

The initial call creates an instance of the class, and defines a number of internal properties, among which one can find the following useful ones:
\begin{itemize}[leftmargin=10pt]
    \item internal copies of the inputs \texttt{XXXeq} and \texttt{DeltaXXX} (and optionally \texttt{angle\_ind}, empty by default)
    \item \texttt{dimension}: problem dimension (number of variables)
    \item \texttt{x}: a YALMIP sdpvar polynomial object, of the dimension of the problem. It represents the variable $\vx$ and is called by the user to define the dynamics of the system
    \item \texttt{t}: sdpvar polynomial of size 1, representing the time variable $t$, which can be needed to define the dynamics of the system (if non-autonomous)
    \item \texttt{solver}, the choice of the solver used in the optimization, defined as Mosek by default, it can also be SeDuMi
    \item \texttt{verbose}, the value of the verbose parameters of the YALMIP optimization calls, defined at 2 by default (all numerical \ac{sdp} solver info displayed), it can also be 1 (selection of info displayed) or 0 (no info about the numerical resolution of the underlying \ac{sdp})
    \item \texttt{dynamics}, a YALMIP polynomial defining the polynomial dynamics $\vf$ of the system.
\end{itemize}

{
Eventually, after declaring all required parameters (see Tab.~\ref{tab:pllparam} and~\ref{tab:smibparam}), the system dynamics are introduced:
\begin{Verbatim}[frame=single, fontsize=\small]
  PLL.dynamics = [PLL.x(2); ...
  -K/tau1*(si(PLL.x(1)) + tau2*co(PLL.x(2))];
\end{Verbatim}
{\color{black}where \texttt{PLL.x} denotes the variable $\vx$ of the PLL system (automatically defined by \texttt{SOStab} when it is first called).} For the sake of brevity the code for SMIB dynamics is reported in Appendix~\ref{sec:smib}. A crucial note is that it features the \textit{recasted} dynamics $\dot\vx = \vf(\vx)$, meaning that $\dot{\theta} = \omega - \omega_s$ is replaced with the dynamics of its recasted images $s = \sin\theta$, $c = \cos\theta$:
\vspace*{-1em}
\begin{subequations}
    \begin{align}
        \dot{s} & = (\omega - \omega_s)c\\
        \dot{c} &= (\omega_s - \omega)s
    \end{align}
\end{subequations}
\vspace*{-1em}
}

\subsection{RoA approximation} \label{sec:prop}
The inner and outer \ac{roa} approximation of the system defined by the call of \texttt{SOStab} are then computed by the methods \texttt{SoS\_in} and \texttt{SoS\_out}, respectively. {For the PLL system, we first compute a degree $16$ outer approximation of} the time $T=1$ \ac{roa} of $\M = \left\{(\phi,\omega) \in \R^2 \; : \; 20\phi^2 + 0.05\omega^2 \leq 1.7^2 \right\}$ using the following commands:
\begin{Verbatim}[frame=single, fontsize=\small]
  T=1;        epsilon = 1.7;          d=16;
  A = [[20^(1/2) 0] ; [0 20^(-1/2)]];
  [vol, vc, wc] = PLL.SoS_out(d,T,epsilon,A);
\end{Verbatim}
{In the SMIB case, we seek a degree $6$ inner approximation of the time $T=10$ \ac{roa} of $\M = \{\vx \in \R^6 \; : \; \|\vx-\vx^\star\| \leq 0.2\}$, hence we enter the commands:
\begin{Verbatim}[frame=single, fontsize=\small]
  T=10;        epsilon = 0.2;          d=6;
  [vol, vc, wc] = SMIB.SoS_in(d,T,epsilon);
\end{Verbatim}
}
Additionnal properties are related to a specific solution of the optimization problem. They are calculated at each call of the optimization and stored until the next call, \textit{i.e.} each of them correponds to the previous optimization call:
\begin{itemize}[leftmargin=10pt]
\item internal copies of the inputs \texttt{d}, \texttt{T}, \texttt{epsilon} (and optional \texttt{A}, set to identity by default); recall that \texttt{d} is the degree of the polynomials $v_d$ and $w_d$
\item \texttt{vcoef\_outer}, the coefficients of the solution $v_d^{out}$ for the last calculated outer approximation of the ROA
\item \texttt{wcoef\_outer}, the coefficients of the solution $w_d^{out}$ for the last calculated outer approximation of the ROA
\item \texttt{vcoef\_inner}, the coefficients of the solution $v_d^{in}$ for the last calculated inner approximation of the ROA
\item \texttt{wcoef\_inner}, the coefficients of the solution $w_d^{in}$ for the last calculated inner approximation of the ROA
\item \texttt{solution}, a volume approximation of the last calculated ROA, \textit{ie} the solution $\lambda_d$ of the optimization problem
\end{itemize}

{
The output of SOStab can be used as follows: as introduced in Section~\ref{sec:roa}, $\lambda_d$ and $w_d$ are used to assess the numerical solver performance; the key output is then the polynomial classifier $v_d$, which is used as follows: in \textit{any} situation where the \texttt{dynamics} are valid for the considered system (be it initialization or post-fault transients) and its state is represented by $\vx(t)$, $v_d^{in}(0,\vx(t)) < 0$ means that the system will remain in the secure zone $\X$ and hit its target: $\vx(t+T) \in \M$; conversely, if $v_d^{out}(0,\vx(t)) < 0$, then either the state will leave $\X$ before time $t+T$ or it will miss its target: $\vx(t+T)\notin \M$; hence, the only remaining uncertainty will be when both $v_d^{in/out}(0,\vx(t)) \geq 0$. Hence, assessing the stability of a configuration $\vx(t)$ boils down to evaluating the polynomial $v_d(0,\vx(t))$, which can be done instantly.
}

\section{Case study}
\label{sec:results}

In this section, we display and comment the results obtained by running the codes presented in \ref{sec:tuto}.

\subsection{RoA of the Phase Locked Loop system}

{When called according to the instruction of Section~\ref{sec:tuto} to assess the stability of the PLL model, SOStab} returns the approximate surface \texttt{vol}$=\lambda_d^{out}$ of the computed \ac{roa} estimate in the phase space, as well as two vectors \texttt{vc} and \texttt{wc} consisting of the coefficients of polynomials $v_d^{out}$ and $w_d^{out}$ respectively. We run SOStab for the PLL system for various values of $d$ and compile the results in Tab.~\ref{tab:d}. As $\lambda_d^{out}$ is a proxy for the size of the outer \ac{roa} estimate, the smaller it is, the more accurate the approximation. Here one can observe that as expected the accuracy increases with $d$ at the price of higher computational time. The CPU times were obtained on a Macbook laptop with an Apple M2 chip and 16 GB of RAM.

\begin{table}[htbp]
    \centering
    \begin{tabular}{|c|c|c|}
        \hline
        $d$ & $\vphantom{\displaystyle\int}\lambda_d^{out}$ & CPU time (s) \\
        \hline
        4 & 4.0000 & 2.8169 \\
        8 & 3.5892 & 4.6729 \\
        12 & 3.1284 & 14.5575 \\
        16 & 2.9346 & 45.2185 \\
        \hline
    \end{tabular}
    \vspace{1em}

    \caption{Outputs of SOStab depending on precision parameter $d$}
    \label{tab:d}
\end{table}

Once the optimization problem is solved, {the following command plots two-dimensional slices of the boundary $$\partial \cR_d^{out} = \{\vx\in\R^n \; : \; v(0,\vx) = 0\}$$} 
\vspace{-1em}
\begin{Verbatim}[frame=single, fontsize=\small]
  PLL.plot_roa(1,2,'outer');
\end{Verbatim}
where the first two arguments indicate the indices of the represented variables (respectively in abscissa and ordinate), in case there are more than two. The string \texttt{'outer'} indicates that the toolbox plots an outer estimate of the \ac{roa}. 

\begin{figure}[htbp]
\centerline{\includegraphics[width=0.7\linewidth, trim=0.1cm 0.1cm 0.1cm 0.65cm,clip]{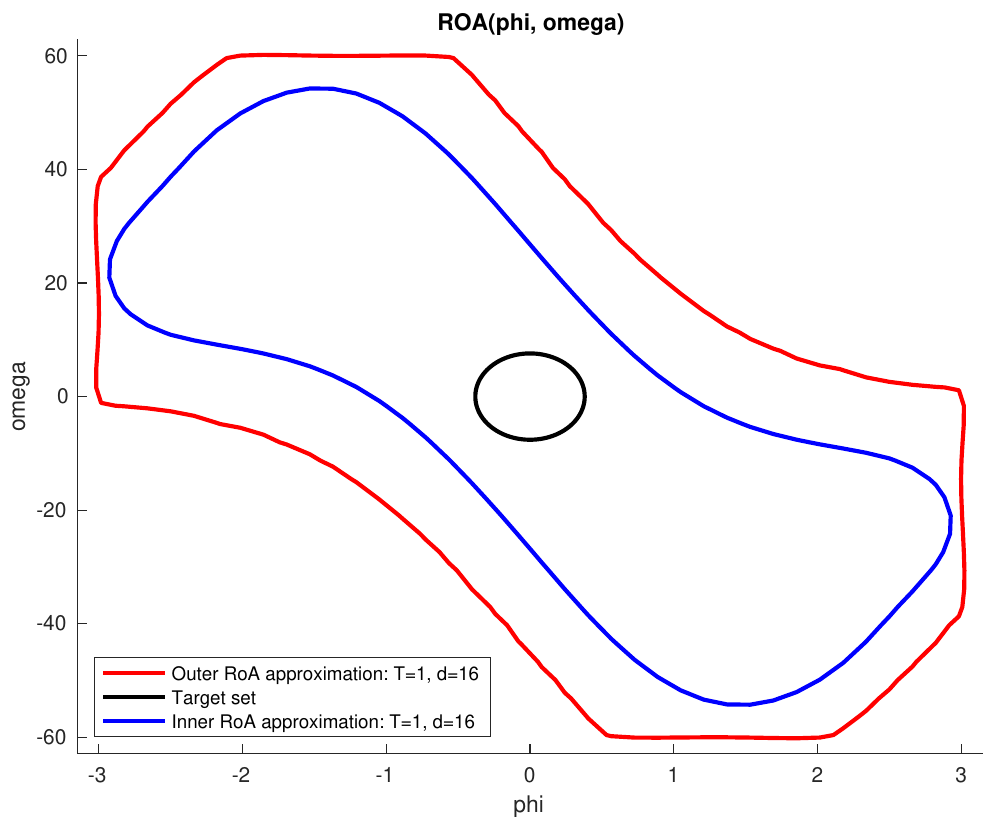}}
\caption{Inner and outer RoA approximations for the PLL system.}
\label{fig:roa-pll}
\end{figure}

For inner \ac{roa} approximation~\cite{kordaInnerROA}, the commands are
\begin{Verbatim}[frame=single, fontsize=\small]
[vol, vc, wc] = PLL.SoS_in(d, T, epsilon);
PLL.plot_roa(1, 2, 'inner',1);
\end{Verbatim}
Here the last argument with value \texttt{1} is an optional argument that asks SOStab to also represent the target set $\M$ in the figure. {The admissible set $\X$ is the whole plotting window, so that it is always visualized.}
This yields the plot represented in Fig. \ref{fig:roa-pll}, which can be compared to~\cite[Fig. 10, Right]{polypll}.

It can also be interesting to display 3D-plots of polynomials $v_d^{out}$ and $w_d^{out}$, which can be performed by the commands (see Fig.~\ref{fig:w-pll}: one retrieves the shape of the inner \ac{roa} estimate):
\begin{Verbatim}[frame=single, fontsize=\small]
  PLL.plot_v(1, 2, 'outer');
  PLL.plot_w(1, 2, 'outer');
\end{Verbatim}
\begin{figure}[htbp]
    \vspace*{-1.5em}
    \centerline{\includegraphics[width=0.8\linewidth, trim=0cm 0cm 0.3cm 1.2cm, clip]{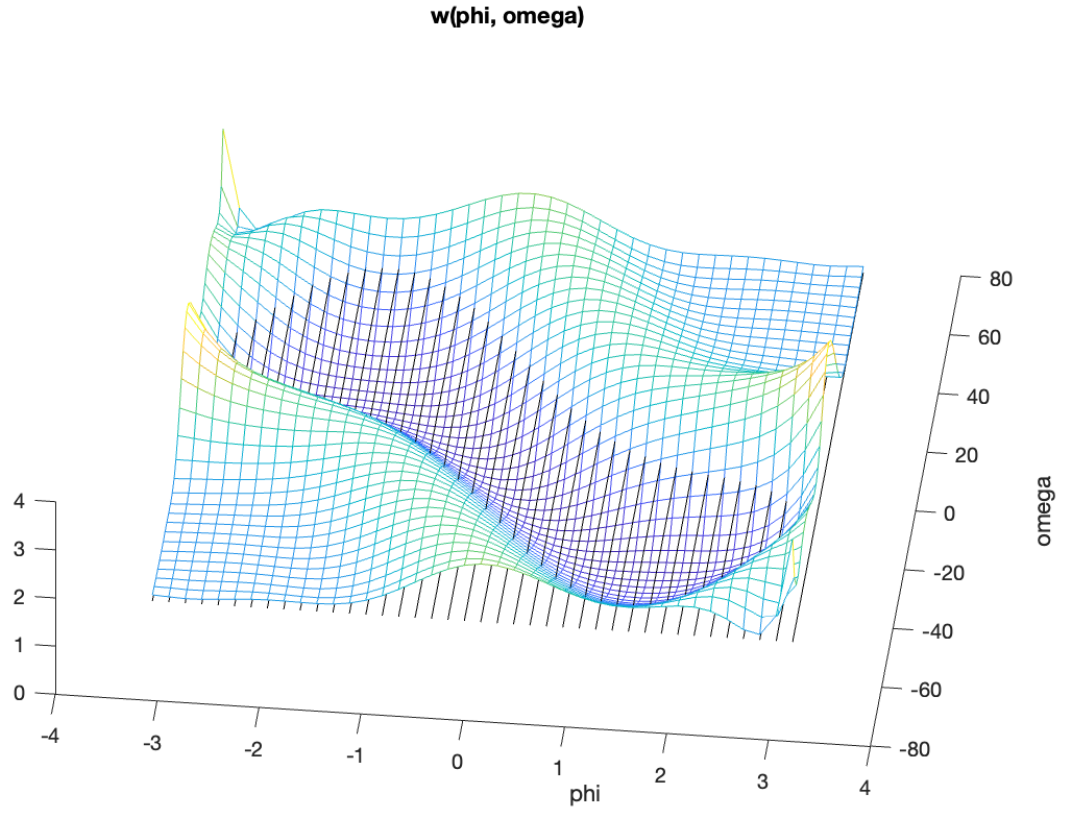}}
    \caption{Plot of $w_d^{in}$ for the PLL system}
    \label{fig:w-pll}
\end{figure}

Of course, one can also represent the certificates $v_d^{in}$ and $w_d^{in}$ obtained in inner approximation, simply by setting the last argument at \texttt{'inner'}. 

{ Those 3D plots can be useful when the validity of the \ac{roa} estimate is unclear. Indeed, although those estimate are theoretically certified, SOStab can fail to compute an \ac{roa} approximation, due to bad conditioning and numerical solver inaccuracy; to detect such behaviour, one can plot the graph of $w_d$: if it is almost flat with values $w_d(\vx) \simeq 1$ everywhere, then the \ac{roa} estimation failed.}

\subsection{RoA of the single machine - infinite bus system}

Instead of performing Taylor expansions as in the previous section, it is also possible to directly tackle trigonometric functions, through the algebraic change of variables in Eq.~\eqref{eq:variable}~\cite{joszTSA}. 

\begin{figure}[htbp]
    \centerline{\includegraphics[width=0.75\linewidth, trim=3cm 1.3cm 4cm 11.3cm,clip]{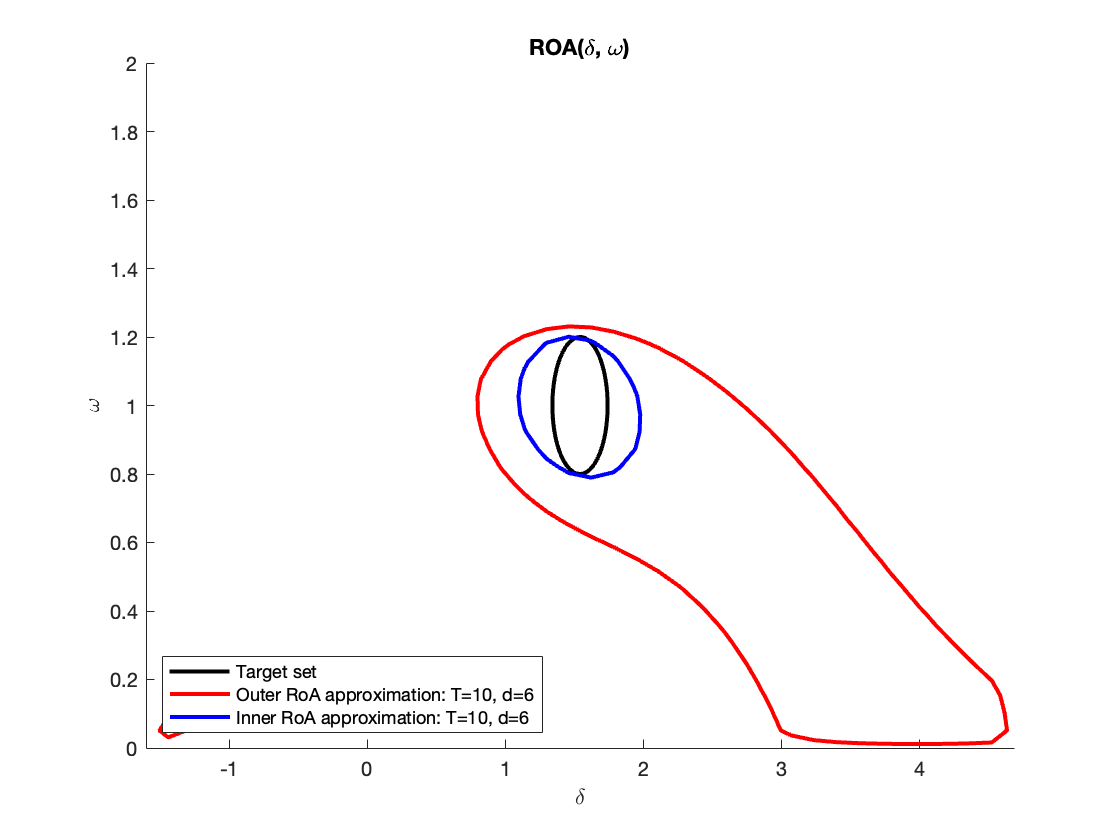}}
    \caption{RoA approximation of the SMIB in the $(\delta,\omega)$ coordinates.}
    \label{fig:smib_freq}
\end{figure}

\begin{figure}[htbp]
    
    \centerline{\includegraphics[width=0.75\linewidth, trim=3.2cm 0.9cm 4.2cm 2.6cm,clip]{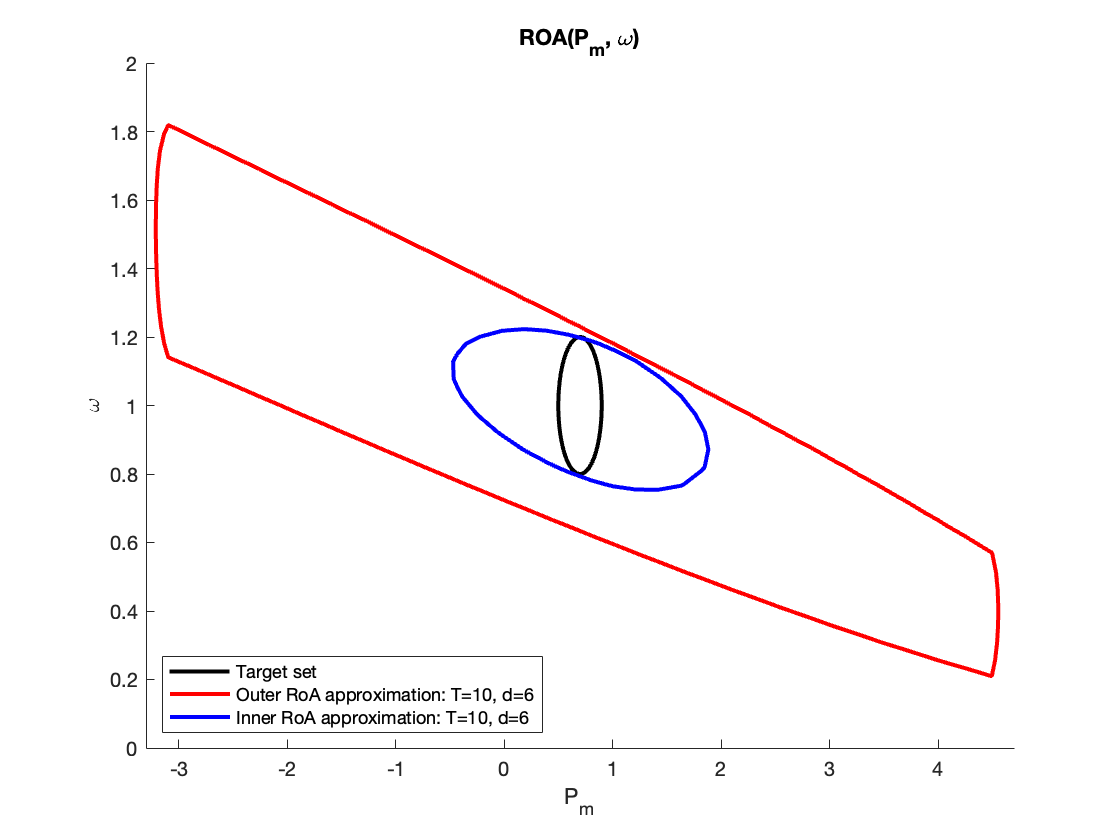}}
    \caption{RoA approximation of the SMIB in the $(P_m,\omega)$ coordinates.}
    \label{fig:smib_pow}
\end{figure}

{After the commands displayed in Section~\ref{sec:tuto}, a user has access to the description of the inner \ac{roa} estimate through the polynomial $v_d^{in}(0,\bfx)$. It is also possible to compute an inner estimate and plot both on a 2D graph in the $(\theta,\omega)$ coordinates (see Fig.~\ref{fig:smib_freq}), with the following commands:
\begin{Verbatim}[frame=single, fontsize=\small]
  SMIB.plot_roa([1,2],3,'inner',1);
  SMIB.SoS_out(d, T, epsilon);
  SMIB.plot_roa([1,2],3,'outer',1);
\end{Verbatim}
Here the argument \texttt{[1,2]} means that the abscissa of the plot should be the phase variable $\theta$, which SOStab knows only from its images $x_1 = \sin \theta$ and $x_2 = \cos \theta$. More generally, setting \texttt{[i,j]} in the coordinate arguments of SOStab results in the abscissa of the plot being the phase variable $\theta$ such that $x_i = \sin \theta$ and $x_j = \cos \theta$. 

The second argument \texttt{3} sets $x_3 = \omega$ as the ordinate of the plot. Another choice such as \texttt{6,3} (without brackets) would result in plot coordinates $(x_6,x_3) = (P_m,\omega)$ as in Fig.~\ref{fig:smib_pow}. Those figures are similar to the existing results~\cite[Fig. 5.b-c]{Marinescu}. For the degree 6 estimates of the SMIB model \ac{roa}, the outer approximation took $49.085$ seconds, while the inner approximation (much more difficult in practice due to conditioning technicalities) took approximately 2 hours.
}

{
\subsection{A curse of dimensionality}
The above examples all include very simplified models of power system devices, with a 2D 2OM PLL and a 5D SMIB model. This is due to the fact that, for computational times and numerical conditioning to be reasonable, the toolbox currently works only with low-dimensional systems. Indeed, SOStab is internally asked to solve LMIs of combinatorial size $\binom{n+d+1}{d}$ where $n$ is the number of variables modelling the system and $d$ is the degree of the certificates $v_d$ and $w_d$. The current computational limits of a standard use of \ac{sosp} in stability analysis were reached in~\cite{Marinescu,oustry} and the 5D SMIB model. A more promising approach for scaling the method consists in splitting the considered problems into subproblems of more modest size, based on the structure of the problem. For instance,~\cite{wang} exploits the algebraic notions of term sparsity and chordal sparsity to handle up to 16 dynamical variables while keeping convergence guarantees. More physics-oriented are the methods in~\cite{tacchi} (resp.~\cite{subotic}), where causality (resp. time scale separation) is exploited to decompose the full problem into more tractable subproblems. Such ideas have been successfully implemented in the simpler (because static) AC-OPF problem setting, for which \textit{continental-scale} systems have been handled by \ac{sosp} methods~\cite{josz,franc}. While very promising, these methodologies pose significant challenges when it comes to their automation into a Matlab toolbox, and hence were not implemented in the current version of SOStab. This, although out of the scope of the present article, will be deeply investigated in future works.
}

\section{Discussions and conclusions}
\label{sec:conclusion}
This work presented a new Matlab Toolbox called SOStab, which aims at helping a non-expert in polynomial optimization to use the frameworks developed in~\cite{henrionOuterROA,kordaInnerROA,joszTSA} through a plug-and-play interface. {Such interface would greatly facilitate the setup of new experiments to leverage the conditioning issues; further, it could be included in future schemes aimed at scaling the method to high dimensions.
}

\subsection{Contribution highlights} \label{sec:function}

{In practice,~\cite{henrionOuterROA,kordaInnerROA,joszTSA} only provide a \textit{methodology} for \ac{roa} approximation. As detailed in those references as well as} Appendix~\ref{sec:hierarchy}, the Moment-\ac{sos} hierarchy consists in solving \ac{sosp} problems of increasing size, which requires to follow a number of steps, related to real algebraic geometry and optimization:
\begin{enumerate}
    \item \label{nb:geo} Defining the geometric characteristics of the problem (polynomial dynamics $\vf$, time horizon $T$, admissible and target sets $\X$ and $\M$)
    \item \label{nb:alg} Defining an algebraic description of these geometric characteristics: polynomials $\vg \in \R[\rmx]^m$ s.t. $\X = \{\vx \; : \; \vg(\vx) \geq 0\}$ and $\vh\in\R[\rmx]^\ell$ s.t. $\M = \{\vx \; : \; \vh(\vx) \geq 0\}$
    \item \label{nb:mom} Coding a method to integrate polynomials over $\X$ (i.e. computing the moments of the Lebesgue measure on $\X$)
    \item \label{nb:sos} Writing the \ac{sosp} problems with explicit positivity constraints; this requires introducing internal \ac{sos} certificates as decision variables
    \item \label{nb:sdp} Recasting \ac{sos} constraints as LMIs, solving the resulting \ac{sdp} problem and converting the solution back to the polynomial framework
    \item \label{nb:rep} Extracting the corresponding certificates $v_d
    ,w_d
    $ and using them to characterize and plot relevant representations of the computed \ac{roa} approximation.
\end{enumerate}
 Existing frameworks~\cite{yalmip,GloptiPoly,SOSTOOLS} have been designed to automate step~\ref{nb:sdp} which appears in all instances of \ac{sosp}. As a result, they 
are very flexible in their use, but they also require their user to perform \mbox{steps~\ref{nb:geo}--\ref{nb:sos}} and~\ref{nb:rep} by hand, which involves solid knowledge in \ac{sosp} and may be time consuming and prone to human errors 
(especially in step~\ref{nb:sos})
; hence their use is usually not smooth, even for experts. 

In contrast, with SOStab, only step~\ref{nb:geo} is left to the user, and all the other operations are automatically performed. 
More precisely, intrinsic properties of the dynamical system are defined as presented in Section~\ref{sec:init}, and then settings for finite horizon \ac{roa} are the input of methods \texttt{SoS\_in} and \texttt{SoS\_out}; 
with this, steps~\ref{nb:alg}--\ref{nb:sdp} are performed through a call to YALMIP, for inner and outer \ac{roa} approximation respectively, and output an optimal value \texttt{solution} and the coefficients of the optimal polynomials $v_d
$ and $w_d
$. 

The benefits of this contribution are the following:
\begin{itemize}[leftmargin=16pt]
    \item[(a)] Knowledge on \ac{sos} programming becomes optional to use the Lasserre hierarchy for \ac{roa} approximation.
    \item[(b)] 
    The user input is significantly reduced, limiting implementation efforts.
    \item[(c)] The toolbox comes with a plug-and-play design that allows one to repeat multiple experiments, reproduce existing results from the literature and solve new problems.
\end{itemize}


\subsection{Limitations and future works}

However, some limitations remain to be leveraged. For instance, while convex, \ac{sdp} problems can be ill-conditioned, which sometimes results in poor numerical behavior with $w^\star_d \simeq 1$ and meaningless plots. It is possible to rescale \ac{sos} constraints to mitigate that phenomenon, although finding the appropriate rescalings is non-trivial. 

Moreover, inner \ac{roa} approximation requires an algebraic representation of the boundary $\partial\X$ of the admissible set $\X$~\cite{kordaInnerROA}, and while choosing a box is the most physically relevant (and the easiest to integrate polynomials on), it induces some numerical difficulties that would not arise if $\X$ were described by a single polynomial. This can be solved either by changing the description of $\partial\X$, or by changing the admissible set $\X$. 


Last but not least, at a given number $n$ of variables (such that $\vx \in \R^n$ and given precision degree $d$, SOStab (and the Lasserre hierarchy in general) requires to solve size $\binom{n+d+1}{d}$ \ac{sdp} problems, which quickly becomes intractable on any computer. To tackle this issue, structure exploiting methods have been developed {in~\cite{tacchiThese,tacchi,subotic,wang,josz,franc}}, which consist in splitting the underlying \ac{sdp} problems into problems of smaller size, while keeping as much computing precision as possible; these techniques bear the potential for scaling the hierarchy up to more realistic power systems such as fully modelled power converters or low order models of distributions network (which exhibit a radial structure one can exploit in computations).

Future works on the SOStab class will include:
\begin{itemize}[leftmargin=16pt]
    \item[(a)] Running the toolbox on more sophisticated case studies such as power converters
    \item[(b)] Improving the inner approximation scheme to increase the accuracy of each relaxation
    \item[(c)] Supporting richer admissible $\X$ and target $\M$ sets, such as ellipsoids, $\ell^p$-balls, annuli...
    \item[(e)] Supporting bases of polynomials other than monomials (Chebyshev, Legendre, trigonometric polynomials)
    \item[(f)] Exporting the toolbox to other softwares compatible with existing \ac{sdp} solvers, such as Julia or Python
    \item[(g)] Adding structure exploiting methods to scale the method to higher dimensional dynamical systems
\end{itemize}

\appendix

\subsection{Lasserre hierarchy for Region of attraction} \label{sec:hierarchy}

In this section, the generic problem of computing the finite time \ac{roa} of a given target set is presented, along with the \ac{sos} framework to address it. Consider the system
\begin{subequations} 
\begin{equation} 
    \dot{\vx} = \vf(\vx)
\end{equation}
with vector field $\vf \in C^\infty(\R^n)^n$ and state constraint
\begin{equation} 
    \vx(t) \in \X
\end{equation}
\end{subequations}
for some subset $\X \subset \R^n$ representing security constraints. 
\begin{dfn} 
Given a time horizon $T \in (0,\infty]$ and a closed target set $\M \subset \X$, the Region of Attraction (RoA) of $\M$ in time $T$ is defined as
\begin{equation} 
    \cR_T^\M := \left\{\vx(0) \in \R^n \; : \begin{array}{c}
        \forall t \in [0,T), \quad \vx(t) \in \X \\
        \dist(\vx(t),\M) \underset{t\to T}{\longrightarrow} 0
    \end{array} \right\}
\end{equation}
\end{dfn}
\begin{rem} Definition~\ref{dfn:RoA} covers many frameworks, such as:
\begin{itemize}[leftmargin=10pt]
    \item Infinite time RoA ($T=\infty$, $\X = \R^n$, often $\M=\{0\}$)
    \item Maximal positively invariant set ($T=\infty$, $\M=\X \varsubsetneq \R^n$)
    \item Constrained finite time RoA ($T < \infty$, $\X$ compact)
\end{itemize}
\end{rem}

We now introduce an infinite dimensional Linear Programming (LP) problem that is related to the constrained finite horizon RoA (see~\cite{henrionOuterROA} for details):
\begin{subequations} \label{eq:LP}
\begin{align}
    W^\star := \inf \; \; & \int_\X w(\vx) \; \rmd\vx \notag \\
    \mathrm{s.t.} \; & v \in C^\infty(\R^{n+1}), w \in C^\infty(\R^n) \notag \\
    & w \geq 0 && \text{on } \X \label{con:pos} \\
    & w \geq v(0,\cdot) + 1 && \text{on } \X \label{con:w_dom_v} \\
    & \cL_{\vf} v := \partial_t v + \vf^\top \vpartial_{\vx} v \leq 0 && \text{on } \bbGamma \label{con:v_decr} \\
    & v(T,\cdot) \geq 0 && \text{on } \M \label{con:target}
\end{align}
\end{subequations}
where $\bbGamma := [0,T]\x\X \Subset \R^{n+1}$ denotes a time-state cylinder.

\begin{prop}[{\cite{henrionOuterROA}}] \label{prop:outer} Let $(v,w)$ be feasible for \eqref{eq:LP}. Then,
\begin{subequations} 
\begin{align}
    \cR_T^\M & \subset \{\vx \in \R^n \; : \; v(0,\vx) \geq 0\} \\ 
    & \subset \{\vx \in \R^n \; : \; w(\vx) \geq 1 \} =: \L(w\geq 1)
\end{align}
\end{subequations}
\end{prop}
With Proposition~\ref{prop:outer}, constraint~\eqref{con:pos} enforces $w \geq \ind_{\cR_T^\M}$, where $\ind_{\A}$ denotes the boolean indicator function of $\A \subset \R^n$ (with value $1$ in $\A$ and $0$ elsewhere). Moreover, it is proven in~\cite{henrionOuterROA} that for any minimizing sequence $(v_d,w_d)_{d\in\N}$ for \eqref{eq:LP}, one has 
$w_d \underset{d\to\infty}{\longrightarrow} \ind_{\cR_T^\M}$ in the sense of $L^1(\X)$, so that the volume of the approximation error $\L(w_d\geq 1) \setminus \cR_T^\M$ converges to $0$.

The Moment-\ac{sos} hierarchy allows its user to compute such a minimizing sequence, under the following assumptions.
\begin{asm} \label{asm:poly} All considered inputs are polynomial:
\begin{enumerate}[label=\ref*{asm:poly}.\arabic*.]
    \item \label{nb:f} $\vf \in \R[\bfx]^n$, so that $v \in \R[\rmt,\bfx] \Longrightarrow \cL_{\vf} v \in \R[\rmt,\bfx]$
    \item \label{nb:X} $\X = \cap_{i=1}^m \L(g_i \geq 0) =: \L(\vg \in \R_+^m)$ with $\vg \in \R[\bfx]^m$
    \item \label{nb:M} $\M = \cap_{j=1}^\ell \L(h_j \geq 0) = \L(\vh \in \R_+^\ell)$ with $\vh \in \R[\bfx]^\ell$
\end{enumerate}
where $\bfx$ (resp. $\rmt$) denotes the dimension $n$ (resp. $1$) indeterminate (i.e. identity function, which can be evaluated in any state $\vx \in \R^n$, resp. time $t \in \R$).
\end{asm}
Then, it is possible to work with polynomial Sums-of-Squares (SoS), with the following definitions.
\begin{dfn} \label{dfn:SoS} Let $p \in \R[\bfx]$. Then, considering $g_0 := 1$,
\begin{itemize}[leftmargin=10pt]
    \item $p$ is SoS \textit{iff} $p = q_1^2 + \ldots + q_N^2$, $q_1,\ldots,q_N \in \R[\bfx]$
    \item $p \in \cQ(\vg)$ \textit{iff} $p = s_0\,g_0 + \ldots + s_m \, g_m$, $s_0,\ldots,s_m$ SoS
    \item $p \in \cQ_d(\vg)$ \textit{iff} $p\in\cQ(\vg)$ with $\max(\deg s_i \, g_i) \leq d$
\end{itemize}
\end{dfn}

Since SoS polynomials are nonnegative by design, it is clear from Definition~\ref{dfn:SoS} that any $p \in \cQ_d(\vg)$ (resp. $\cQ_d(\vh)$) is nonnegative on $\X$ (resp. $\M$), which gives access to a strenghtening of problem \eqref{eq:LP}:
\begin{subequations} \label{eq:SoS}
\begin{align}
    W^\star_d := \inf \; \; & \int_\X w(\vx) \; \rmd\vx \notag \\
    \mathrm{s.t.} \; & v \in \R[\rmt,\bfx], w \in \R[\bfx] \hspace*{-3em} \notag \\
    & w && \in \cQ_d(\vg) \label{con:pos_d} \\
    & w - v(0,\bfx) - 1 && \in \cQ_d(\vg) \label{con:w_dom_v_d} \\
    & -\cL_{\vf} v && \in \cQ_d(\vg, (T-\rmt)\,\rmt) \label{con:v_decr_d} \\
    & v(T,\bfx) && \in \cQ_d(\vh) \label{con:target_d}
\end{align}
\end{subequations}
where $\bbGamma = \L((T-\rmt)\,\rmt \geq 0) \x \X = \L((\vg, (T-\rmt)\,\rmt) \in \R_+^{m+1})$. 

\noindent Problem \eqref{eq:SoS} consists in looking for feasible $(v,w)$ for \eqref{eq:LP} under the form of polynomials, restricting inequality constraint (\ref{eq:LP}x) into SoS constraint (\ref{eq:SoS}x), x$=$a--d. The advantage of this new problem is that the decision variables are now finite dimensional vectors of coefficients, and the SoS constraints can be recast as LMIs~\cite{LasserreBook}. Thus, assuming knowledge of the moments of the Lebesgue measure on $\X$ (i.e. being able to integrate polynomials on $\X$, e.g. if $\X$ is a ball or a box), one is able to solve \eqref{eq:SoS} on a standard computer, provided that $n$ and $d$ are small enough to make the LMIs tractable.

As the new feasible space is strictly included in the former, there is a relaxation gap: $\forall d \in 2\N, \ W^\star < W^\star_d$. However, it is proved in~\cite{henrionOuterROA} that, if $\X$ and $\M$ are bounded, $W^\star_d \underset{d\to\infty}{\longrightarrow} W^\star$. 
Thus, solving instances of \eqref{eq:SoS} gives access to converging outer approximations of the RoA.

\subsection{A standard polynomial system} \label{sec:pol}
To further illustrate how SOStab is able to reproduce existing results in \ac{sosp} for stability analysis, one can consider a reversed-time Van der Pol oscillator, as in \cite{kordaInnerROA}:
\begin{equation}
\begin{pmatrix}
\dot{x_1} \\ \dot{x_2}
\end{pmatrix} = \begin{pmatrix}
-2x_2 \\ 0.8x_1 + 10(x_1^2 - 0.21)x_2
\end{pmatrix}
\end{equation}

The dynamics are polynomial and the stable equilibrium $\vx^\star$ of the system is at the origin, so that it takes very few lines of code to get interesting results:

\begin{Verbatim}[frame=single, fontsize=\small]
VdP = SOStab([0;0], [1.1;1.1]);
VdP.dynamics= [-2*VdP.x(2); 0.8*VdP.x(1)
+ 10*(VdP.x(1)^2 - 0.21)*VdP.x(2)];
d = 12;         T = 1;         epsilon = 0.5;
\end{Verbatim}

\begin{Verbatim}[frame=single, fontsize=\small]
[vol, vc, wc] = VdP.SoS_out(d, T, epsilon);
VdP.plot_roa(1, 2, 'outer');
[vol, vc, wc] = VdP.SoS_in(d, T, epsilon);
VdP.plot_roa(1, 2, 'inner',1);
\end{Verbatim}

\begin{figure}[htbp]
\centerline{\includegraphics[width=0.8\linewidth, trim=1.8cm 7.9cm 3.6cm 9.5cm, clip]{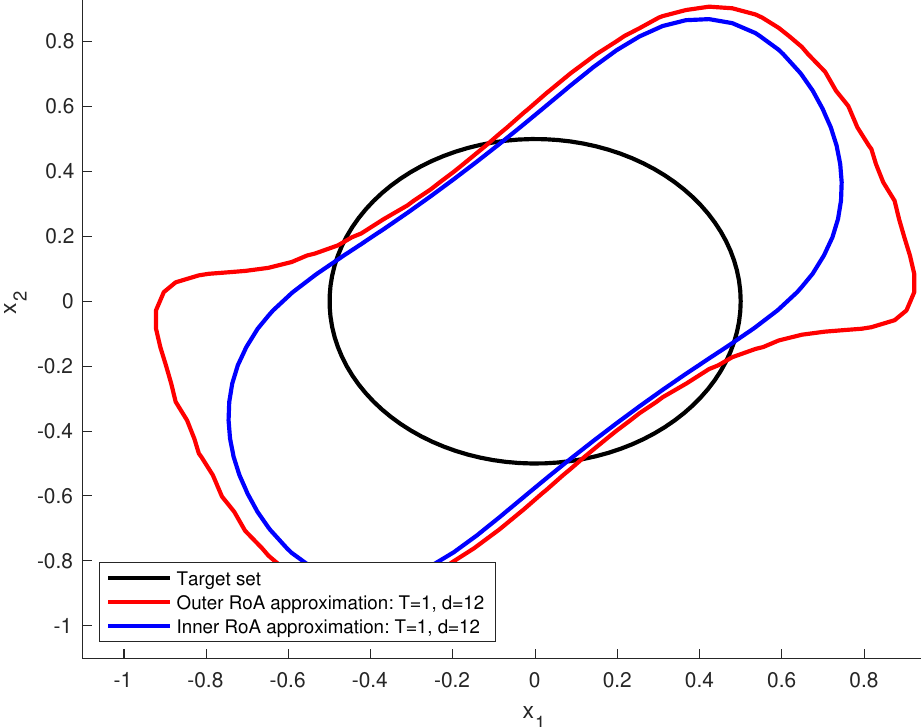}}
\caption{Inner and outer RoA approximations for the Vanderpol system.}
\label{fig:roa-vdp}
\end{figure}

This gives the plot represented in Fig. \ref{fig:roa-vdp}, which reproduces results presented in~\cite[Figure 2]{henrionOuterROA} and~\cite[Figure 3]{kordaInnerROA}.

\subsection{Code for the SMIB RoA dynamics}
\label{sec:smib}
\begin{Verbatim}[frame=single, fontsize=\small]
  Td = 9.67;  xd = 2.38;  xpd = 0.336;
  xq = 1.21;  H = 3;  r = 0.002;  w = 1;
  R = 0.01;  X = 1.185;  V = 1;  Ta = 1;
  Ka = 70;  Tg = 0.4;  Kg = 0.5;  P = 0.7;
  
  iq = ((X+xpd)*V*SMIB.x(1) - (R+r)*...
       (V*SMIB.x(2)-SMIB.x(4)))/...
       ((R+r)^2 + (X+xpd)*(X+xq));
  id = (X+xq)/(R+r)*iq - 1/(R+r)*V*SMIB.x(1);
  vd = xq*iq - r*id;
  vq = R*iq + X*id + V*SMIB.x(2);
  Vt = vd^2+vq^2;
  
  SMIB.dynamics = [(SMIB.x(3)-weq)*...
     SMIB.x(2); (weq-SMIB.x(3))*SMIB.x(1);...
     (SMIB.x(6) - vd*id - vq*iq - r*id^2...
       - r*iq^2)/(2*H);...
  (SMIB.x(5) - SMIB.x(4) + (xpd-xd)*id)/Td;...
  (Ka*(V-Vt) - SMIB.x(5))/Ta ;...
  (P - SMIB.x(6) + Kg*(w-SMIB.x(3)))];
\end{Verbatim}

\end{document}